\documentclass[french,,final]{svjour2}

\usepackage[english]{babel}
\usepackage{latexsym}

\usepackage{amsmath}
\usepackage{amsfonts}
\usepackage{amssymb}
\usepackage{verbatim}
\usepackage{psfrag}
\usepackage{graphicx}
\newcommand{\proofend}{\hfill $\square$}
\newcommand{\Z}[1]{\mathbb{Z}/#1\mathbb{Z}}
 
\newcommand{\defn}[1]{{\em #1}}

\newcommand{\K}{\mathbb{C}}

\newcommand{\PGL}{\mathrm{PGL}}

\newcommand{\GL}{\mathrm{GL}}
\newcommand{\im}{\mathbf{i}}
\newcommand{\Pn}{\mathbb{P}^2}
\newcommand{\Aut}{\mathrm{Aut}}
\newcommand{\Sym}{\mathrm{Sym}}

\newcommand{\titreProp}[1]{{\bf - #1} }
\newcommand{\rkPic}[1]{\mathrm{rk\ Pic}(#1)}
\newcommand{\Pic}[1]{\mathrm{Pic}(#1)}
\newcommand{\vspm}{\vspace{-0.05cm}}
 \newcommand{\sfrac}[2]{\leavevmode\kern.1em
            \raise.5ex\hbox{\footnotesize #1}\kern-.1em
                    /\kern-.15em\lower.25ex\hbox{\footnotesize #2}}

\title{{\bf {\Large The number of conjugacy classes of elements of the Cremona group of some given finite order} }}
\author{J\'er\'emy Blanc}
\institute{J\'er\'emy Blanc \at 
Laboratoire J.A. Dieudonn\'e (UMR 6621),
Universit\'e de Nice Sophia Antipolis - C.N.R.S., 
Facult\'e des Sciences - Parc Valrose,
06108 Nice cedex 2 (France),\\ \email{blancj@unice.fr}}

\begin{document}
\maketitle
\begin{abstract}This note presents the study of the conjugacy classes of elements of some given finite order $n$ in the Cremona group of the plane. In particular, it is shown that the number of conjugacy classes is infinite if $n$ is even, $n=3$ or $n=5$, and that it is equal to $3$ (respectively $9$) if $n=9$ (respectively $15$), and is exactly $1$ for all remaining odd orders.\\ Some precise representative elements of the classes are given.\end{abstract}

\section{Introduction}
Let us recall that a \emph{rational transformation} of $\mathbb{P}^2(\K)$ is a map of the form $$(x:y:z) \dasharrow (\varphi_1(x,y,z):\varphi_2(x,y,z):\varphi_3(x,y,z)),$$ where $\varphi_1,\varphi_2,\varphi_3 \in \K[x,y,z]$ are homogenous polynomials of the same degree. If such a map has an inverse of the same type, we say that it is \emph{birational}.

The \defn{Cremona group} is the group of birational transformations of $\mathbb{P}^2(\K)$.
This group has been studied since the $\mathrm{XIX}^{\mathit{th}}$ century by many mathematicians. One of the first natural questions that we may ask when we study some group is the following:
\begin{question}\label{TheQuestion}Given some positive integer $n$, how many conjugacy classes of elements of order $n$ exist in the Cremona group?\end{question}

First of all, it is important to note that the number of conjugacy classes is at least one, for any integer $n$, as the linear automorphism $$(x:y:z)\mapsto (x:y:e^{2\im\pi/n}z)$$ is a representative element of one class. It was proved in \cite{bib:BeB} that all the linear automorphisms of the plane of the same finite order are birationally conjugate (the same is true in any dimension, see \cite{bib:JB}, Proposition 5); to find more conjugacy classes we have therefore to show the existence of non-linearizable birational transformations.

The first answer to Question \ref{TheQuestion} was given in \cite{bib:EB} for $n=2$. Infinitely many involutions which are not conjugate are found. Since the proof of \cite{bib:EB} is considered as incomplete, a precise and complete one may be found in \cite{bib:BaB}.

In \cite{bib:deF}, the answer for $n$ prime is given. It is shown that the number of conjugacy classes is infinite for $n=3,5$ and is equal to $1$ if $n$ is a prime integer $\geq 7$.

For other orders, a lot of examples have been given in the ancient articles (for example in \cite{bib:SK}, \cite{bib:Wim}) and in many more recent articles, the most recent one being \cite{bib:Dol}. However, the precise answer to Question \ref{TheQuestion} was not given for $n$ not prime.

In this paper, we answer to Question \ref{TheQuestion} for any integer $n$, proving the following theorems:
\begin{theorem}
\label{Thm:InfiniteEven35}
For any even integer $n$, the number of conjugacy classes of elements of order $n$ in the Cremona group is infinite. This is also true for $n=3,5$.
\end{theorem}

\begin{theorem}
\label{Thm:FiniteOdd}
For any odd integer $n\not=3,5$, the number of conjugacy classes of elements of order $n$ in the Cremona group is finite.\\
Furthermore this number is equal to 3 (respectively 9) if $n=9$ (respectively if $n=15$) and is $1$ otherwise.
\end{theorem}

\begin{remark}We already stated some of these results in \cite{bib:JBTh}, Theorem 1.\end{remark}

\section{Automorphisms of rational surfaces}
\label{Sec:AutRat}
Let us remark the obvious but important observation: take some birational transformation $\varphi$ of a rational surface $S$. Any birational map $\lambda: S\dasharrow \Pn$ conjugates $\varphi$ to the birational transformation $\varphi_{\lambda}= \lambda\circ\varphi\circ\lambda^{-1}$ of $\Pn$. Although $\varphi_{\lambda}$ is not unique, all the possible $\varphi_{\lambda}$'s form an unique conjugacy class of birational transformations of $\Pn$.

Conversely, taking some birational transformation of $\Pn$, we may conjugate it to a birational transformation of any rational surface. If the order of the transformation is finite, we may furthermore conjugate it to a \emph{(biregular) automorphism} of a rational surface. (See for example \cite{bib:DFE}, Theorem 1.4).

An important family of rational surfaces are the rational surfaces with an ample anticanonical divisor, i.e.\ the
\emph{Del Pezzo surfaces}. These surfaces are $\mathbb{P}^1\times\mathbb{P}^1$, $\mathbb{P}^2$, and the blow-up of $1\leq r\leq 8$ points of $\mathbb{P}^2$ in a general position (i.e.\ such that no irreducible curve of self-intersection $-2$ belongs to the surface).  There is an extensive literature about this; some descriptions may be found for example in \cite{bib:Kol}. Note that the \emph{degree} of such a surface is the square of its canonical divisor, and is an integer between $1$ and $9$; it is $9$ for $\Pn$, $8$ for $\mathbb{P}^1\times\mathbb{P}^1$ and $9-r$ for the blow-up of $r$ points in $\Pn$. Almost all of our examples of rational surfaces will be Del Pezzo surfaces.

\section{Elements of order $3$,$5$ and of any even order - The proof of Theorem \ref{Thm:InfiniteEven35}}
Let us give families of conjugacy classes of elements of order $2$, $3$ and $5$ of the Cremona group.
\begin{example}{\it - Birational transformations of order $2$}\\
\label{Exa:DeJinv}
Let $a_1,...,a_n,b_1,...,b_n \in \K$ be all distinct.  The birational map
$$\Big((x_1:x_2),(y_1:y_2)\Big) \dasharrow \Big((x_1:x_2),(y_2 \prod_{i=1}^n (x_1-b_ix_2):y_1 \prod_{i=1}^n(x_1-a_i x_2))\Big)$$
of $\mathbb{P}^1\times\mathbb{P}^1$ is an involution, which is classically called \emph{de Jonqui\`eres involution}. Its fixed points form a smooth curve $\Gamma\subset \mathbb{P}^1\times\mathbb{P}^1$ of equation $$(y_1)^2\cdot\prod_{i=1}^n(x_1-a_i x_2)=(y_2)^2\cdot\prod_{i=1}^n (x_1-b_ix_2).$$
The restriction to $\Gamma$ of the projection $\mathbb{P}^1\times\mathbb{P}^1\rightarrow \mathbb{P}^1$ on the first factor is a surjective morphism $\Gamma \rightarrow \mathbb{P}^1$ of degree $2$, ramified over the points $(a_1:1),...,(a_n:1),(b_1:1),...,(b_n:1)$. The curve $\Gamma$ is therefore an \defn{hyperelliptic curve}. \\
{\it 
These involutions are birationally equivalent to those of \cite{bib:BaB}, Example~2.4(c).}
\end{example}

\bigskip

\begin{example}{\it - Birational transformations of order $3$}\\ 
\label{Exa:Order3}
Let $F$ be a non-singular form of degree $3$ in $3$ variables and let $\Gamma=\{(x:y:z) \in \mathbb{P}^2\ |\ F(x,y,z)=0\} $ be the smooth cubic plane curve associated to it.

The surface $S=\{(w:x:y:z) \in \mathbb{P}^3\ |\ w^3=F(x,y,z)\} \subset \mathbb{P}^3$ is thus a smooth cubic surface in $\mathbb{P}^3$, which is rational (it is a Del Pezzo surface of degree $3$, see for example \cite{bib:Kol}, Theorem III.3.5). The map $w \mapsto e^{2\im\pi/3}w$ gives rise to an automorphism of $S$ whose set of fixed points is isomorphic to the elliptic curve $\Gamma$.\\
{\it 
Such elements generate cyclic groups of order $3$, already given in  \cite{bib:deF}, Theorem A, case A1.}
\end{example}

\bigskip

\begin{example}{\it - Birational transformations of order $5$}\\
\label{Exa:Order5}
Let us choose $\lambda,\mu \in \K$ such that the surface $$S=\{(w:x:y:z) \in \mathbb{P}(3,1,1,2) \ |\ w^2=z^3+\lambda x^4z+x(\mu x^5+y^5)\}$$ is smooth. The surface $S$ is thus rational (it is a Del Pezzo surface of degree $1$, see \cite{bib:Kol}, Theorem III.3.5) and the map $y\mapsto e^{2\im\pi/5}y$ gives rise to an automorphism of $S$ whose set of fixed points is the union of the point $(1:0:0:1)$ and the trace of the equation $y=0$ on $S$, which is an elliptic curve.\\
{\it 
The corresponding cyclic groups of order $5$ were given in \cite{bib:deF}, Theorem A, case A3.}
\end{example}

To prove Theorem \ref{Thm:InfiniteEven35}, it remains to show the existence of infinitely many conjugacy classes of elements of order $n$, for any even integer $n \geq 4$.
These elements are roots of de Jonqui\`eres involutions (Example \ref{Exa:DeJinv}) and belong to the classical \defn{de Jonqui\`eres group}, which is a subgroup of the Cremona group. We now introduce this group.

\begin{example}{\it - The de Jonqui\`eres group}\\
\label{Exa:DeJg}
The de Jonqui\`eres group is isomorphic to $\PGL(2,\K(x))\rtimes \PGL(2,\K)$, where $\PGL(2,\K)$ acts naturally on $\K(x)$, as  $\PGL(2,\K)=\Aut(\mathbb{P}^1)$ is the automorphism group of $\mathbb{P}^1$ and $\K(x)=\K(\mathbb{P}^1)$ is its function field. To the element $${\Big(}\left(\begin{array}{cc}\alpha(x) & \beta(x)\\ \gamma(x) & \delta(x)\end{array}\right),\left(\begin{array}{cc}a & b\\ c & d\end{array}\right)\Big)\in \PGL(2,\K(x))\rtimes \PGL(2,\K)$$ we associate the birational map $$(x,y) \dasharrow \Big(\frac{ax+b}{cx+d},\frac{\alpha(x)y+\beta(x)}{\gamma(x)y+\delta(y)}\Big)$$
of $\K^2$. The natural inclusion $\K^2 \subset \mathbb{P}^2(\K)$ (respectively $\K^2 \subset \mathbb{P}^1(\K)\times\mathbb{P}^1(\K)$) sends the de Jonqui\`eres group on the group of birational transformations of $\mathbb{P}^2$ (respectively of $\mathbb{P}^1\times\mathbb{P}^1$) that leave invariant the pencil of lines of $\mathbb{P}^2$ passing through one point (respectively that leave invariant one of the two standard pencils of lines of $\mathbb{P}^1\times\mathbb{P}^1$).
\end{example}

In this context, we may look at the subgroup of the de Jonqui\`eres group that fixes some hyperelliptic curve:
\begin{example}{\it - The group of birational transformations that fix some curve}\\
\label{Exa:Jg}
Let $g(x) \in \K(x)^{*}$ be some element which is not a square in $\K(x)$. We denote by $\mathcal{J}_g$ the torus of $\PGL(2,\K(x))$ which is the image in $\PGL(2,\K(x))$ of the subgroup 
$$\mathcal{T}_g=\Big\{\left(\begin{array}{cc}\alpha(x) & \beta(x)g(x)\\ \beta(x) & \alpha(x)\end{array}\right)\ \Big|\ \alpha(x),\beta(x) \in \K(x), \alpha\not=0 \mbox{ or }\beta\not=0\Big\}$$
 of $\GL(2,\K(t))$. The group $\mathcal{J}_g$ corresponds to the group of birational transformations of the form $$(x,y) \dasharrow \Big(x,\frac{\alpha(x)y+\beta(x)g(x)}{\beta(x)y+\alpha(x)}\Big).$$ Note that $\mathcal{T}_g$ is isomorphic to the multiplicative group of the field $$\K(x)\big[\sqrt{g(x)}\big]=\big\{\alpha(x)+\beta(x)\sqrt{g(x)}\ \big|\ \alpha(x),\beta(x) \in \K(x)\big\}.$$ In the case where $g(x)$ is a polynomial without multiple roots, the field $\K(x)[\sqrt{g(x)}]$ is the function field $\K(\Gamma)$ of the smooth curve $\Gamma$ of equation $y^2=g(x)$, and the group $\mathcal{J}_g=\K(\Gamma)^{*}/\K^{*}$ is the group of birational maps of the de Jonqui\`eres group that fix the curve $\Gamma$. (If the degree of $g(x)$ is at least $5$, it is in fact the group of birational maps that fix the curve.)
\end{example}

\bigskip

\begin{proposition}
\label{Prp:ExistenceRoots}
Let $n\geq 1$ be some integer, and let $g(x) \in \K(x)^{*}$ be such that $g(e^{2\im\pi/n}\cdot x)=g(x)$. There exists $\nu(x)\in\K(x)$ such that the $n$-th power of the birational map
 $$\varphi:(x,y) \dasharrow \Big(e^{2\im\pi/n}\cdot x,\frac{\nu(x)y+g(x)}{y+\nu(x)}\Big)$$ is the de Jonqui\`eres involution
$$\varphi^n:(x,y) \dasharrow \big(x,\frac{g(x)}{y}\big).$$
\end{proposition}
\begin{proof}
Note that choosing any $\nu(x) \in \K(x)$, the associated map $\varphi$ belongs to the de Jonqui\`eres group (Example \ref{Exa:DeJg}) and is the composition of $(x,y) \mapsto (e^{2\im\pi/n}\cdot x,y)$ with an element of $\mathcal{J}_g$ defined in Example \ref{Exa:Jg}. The $n$-th power of $\varphi$ in the de Jonqui\`eres group $\PGL(2,\K(t))\rtimes \PGL(2,\K)$ is equal to
$$\Big(\left(\begin{array}{cc}\nu(x) & g(x)\\ 1 & \nu(x)\end{array}\right)
\left(\begin{array}{cc}\nu(\xi\cdot x) & g(x)\\ 1 & \nu(\xi\cdot x)\end{array}\right)
\cdots\left(\begin{array}{cc}\nu(\xi^{n-1}\cdot x) & g(x)\\ 1 &\nu(\xi^{n-1}\cdot x)\end{array}\right),\left(\begin{array}{cc}1 & 0\\ 0 &1\end{array}\right)\Big),$$
where $\xi=e^{2\im\pi/n}$. 
Since the element $\left(\begin{array}{cc}\nu(x) & g(x)\\ 1 & \nu(x)\end{array}\right)\in\mathcal{J}_g\subset\PGL(2,\K(t))$ is the image of some element of $\mathcal{T}_g \subset \GL(2,\K(t))$ corresponding to 
$$\zeta=(\nu(x)+\sqrt{g(x)})\in\K(x)[\sqrt{g(x)}]^{*},$$
the element $\varphi^n\in\mathcal{J}_g\subset\PGL(2,\K(t))$ is therefore the image of the element 
$$\zeta\cdot \sigma(\zeta)\cdot\sigma^2(\zeta)\cdots\sigma^{n-1}(\zeta)\in \K(x)[\sqrt{g(x)}]^{*},$$
where $\sigma$  is the automorphism of $\K(x)[\sqrt{g(x)}]^{*}$ that sends $x$ on $\xi x$ and acts trivially on $\K[\sqrt{g(x)}]$. Let us look at the morphism 
$$\begin{array}{rccc}N:&\K(x)[\sqrt{g(x)}]^{*}&\longrightarrow& \K(x)[\sqrt{g(x)}]^{*}\\ & \tau &\mapsto& \tau\cdot \sigma(\tau)\cdot \sigma^2(\tau)\cdots \sigma^{n-1}(\tau).\end{array}$$
All elements of its image are invariant by $\sigma$ and thus belong to the multiplicative group of the field  $\K(x^n)[\sqrt{g(x)}]$. Furthermore, the map $N$ is the norm of the field extension $\K(x)[\sqrt{g(x)}]/{\K(x^n)[\sqrt{g(x)}]}$.   

Since this is a finite Galois extension, and the field $\K(x)$ has the $\mathrm{C}_1$-property (by Tsen theorem), the norm $N$ is surjective (see \cite{bib:Serre}, X.7, Propositions 10 and 11). We may thus choose an element $\zeta_0=\alpha(x)+\beta(x)\sqrt{g(x)}$ whose norm is equal to $\sqrt{g(x)}$. As $\beta(x)$ is certainly not equal to zero, we may choose $\nu(x)=\alpha(x)/\beta(x)$, so that $\zeta=\zeta_0/\beta(x)$ is sent by $N$ on $\beta\cdot \sqrt{(g(x)}$, whose image in $\PGL(2,\K(t))$ is $\left(\begin{array}{cc}0 & g(x)\\ 1 &0\end{array}\right)$, as we wanted.\proofend
\end{proof}

\bigskip
We give now explicitly a family of examples produced in Proposition \ref{Prp:ExistenceRoots}.
\begin{example}
\label{Exa:Root2oddorder}
Let $n=2m$, where $m$ is an odd integer and let $h \in \K(x)$ be a rational function. We choose $\alpha$ to be the birational transformation
$$\alpha:(x,y)\dasharrow \Big(e^{2\im\pi/n}\cdot x, \frac{h(x^m)y-h(x^m)h(-x^m)}{y+h(x^m)}\Big).$$
Compute $\alpha^2:(x,y)\dasharrow (e^{2\im\pi/m}\cdot x, \frac{-h(x^m)\cdot h(-x^m)}{y})$ and see that this is the composition of the commutative birational transformations $(x,y) \mapsto (e^{2\im\pi/m}\cdot x,y)$ and $(x,y)\dasharrow (x, \frac{-h(x^m)\cdot h(-x^m)}{y})$ of order respectively $m$ and $2$. Thus, the order of $\alpha^2$ is $2m=n$ and $\alpha^{n}=\alpha^{2m}$ is the birational involution
$$\alpha^{n}:(x,y)\dasharrow \Big(x, \frac{-h(x^m)\cdot h(-x^m)}{y}\Big).$$ 
\end{example}

\bigskip

We are now able to prove Theorem  \ref{Thm:InfiniteEven35}, i.e.\ to show the existence of infinitely many conjugacy classes of elements of order $n$ in the Cremona group, for any even integer $n$ and for $n=3,5$.

\begin{proof}{\it of Theorem \ref{Thm:InfiniteEven35}}\\
First of all, note that a birational transformation sends a non-rational curve on another non-rational curve, which have an isomorphic normalisation (the same result for rational curve is false, as these curves may be collapsed on one point). If two birational transformations $\alpha$, $\beta$ are conjugate by $\varphi$, the element $\varphi$ sends the non-rational curves fixed by $\alpha$ on the non-rational curves fixed by $\beta$. (In fact there is at most one such curve, but we will not need it here).

Choosing different de Jonqui\`eres involutions (Example \ref{Exa:DeJinv}), the possible curves fixed are all the hyperelliptic curves. As the number of isomorphism classes of such curves is infinite, we obtain infinitely many conjugacy classes of de Jonqui\`eres involutions in the Cremona group. (In fact, there exist some other families, called Geiser and Bertini involution, see \cite{bib:BeB}).

The same arguments works for elements of order $3$ and $5$ (Examples \ref{Exa:Order3} and \ref{Exa:Order5}), that may fix all the elliptic curves, whose number of isomorphism classes is also infinite.

Taking $n\geq 2$, and any polynomial $g \in \K[x^n]$, there exists an element $\alpha$ in the Cremona group which has order $2n$ and such that $\alpha^n$ is the birational involution $(x,y) \dasharrow (x,\frac{g(x)}{y})$ (Proposition \ref{Prp:ExistenceRoots}). As this involution fixes the hyperelliptic curve $y^2=g(x)$, the number of conjugacy classes of such elements (when changing the element $g$) is infinite.\proofend
\end{proof}

\section{Elements of odd order $\geq 7$ - The proof of Theorem \ref{Thm:FiniteOdd}}
As it was said in Section \ref{Sec:AutRat}, any birational transformation of finite order of the plane is conjugate to an  automorphism $g$ of some rational surface $S$. We may then assume that the pair $(g,S)$ is minimal and use the following result, proved in \cite{bib:Man}.
\begin{proposition}\label{Prp:Man}
Let $g$ be some automorphism of a rational surface $S$, such that the pair $(g,S)$ is minimal  (i.e.\ every $g$-equivariant birational morphism $S\rightarrow S'$ is an isomorphism). Then, one of the two following cases occurs:
\begin{itemize}
\item
$\rkPic{S}^{g}=1$ and $S$ is a Del Pezzo surface.
\item
$\rkPic{S}^{g}=2$ and $g$ preserves a conic bundle structure on $S$.\\ (i.e.\ there exists some morphism $\pi:S\rightarrow \mathbb{P}^1$ with fibres isomorphic to $\mathbb{P}^1$, except for a finite number of singular curves, that consist on the union of two intersecting curves isomorphic to $\mathbb{P}^1$; and $g$ sends any fibre of $\pi$ on another fibre). \end{itemize}
\end{proposition}

\bigskip

To prove Theorem \ref{Thm:FiniteOdd}, we enumerate the possibilites of pairs $(g,S)$ where $g$ is an automorphism of odd order $\geq 7$, using Proposition \ref{Prp:Man}.
The following lemma will help us to prove the Theorem for Del Pezzo surfaces:
\begin{lemma} \titreProp{Size of the orbits}\\
\label{Lem:SizeOrbits}
Let $S$ be a Del Pezzo surface, which is the blow-up of $1\leq r\leq 8$ points of\hspace{0.2 cm}$\mathbb{P}^2$ in general position,
and let $G \subset \Aut(S)$ be a finite subgroup of automorphisms with $\rkPic{S}^{G}=1$. Then:
\begin{itemize}
\item
$G\not=\{1\}$;
\item
the size of any orbit of the action of $G$ 
on the set of exceptional divisors is divisible by the degree of $S$, which is $(K_S)^2$;
\item
in particular, the order of $G$ is divisible by the degree of $S$.
\end{itemize}
\end{lemma}\begin{proof}
It is clear that $G \not=\{1\}$, since $\rkPic{S}>1$.
Let $D_1,D_2,...,D_k$ be $k$ exceptional divisors of $S$, forming an orbit of $G$. 
The divisor $\sum_{i=1}^{k} D_i$ is fixed by $G$ and thus is a multiple of $K_S$. 
We can write $\sum_{i=1}^{k} D_i=a K_S$, for some rational number $a \in \mathbb{Q}$. In fact, since $a K_S$ is effective, we have $a<0$; furthermore $a \in \mathbb{Z}$, since the canonical divisor is not a multiple in $\Pic{S}$.
The $D_i$'s being irreducible and rational, we deduce from the adjunction formula $D_i (K_S +D_i)=-2$ that $D_i\cdot K_S=-1$.
Hence \begin{center}$K_S \cdot \sum_{i=1}^{k} D_i=\sum_{i=1}^{k}K_S \cdot D_i=-k=K_S \cdot a K_S=a (K_S)^2$.\end{center} Consequently, the degree of $S$ divides the size $k$ of the orbit. \proofend
\end{proof}

\bigskip

We decompose now our investigations on different surfaces.

\begin{proposition}
\label{Prp:AutP1P1conjP2}
Any automorphism of $\mathbb{P}^1\times\mathbb{P}^1$ is birationally conjugate to a linear automorphism of $\Pn$.
\end{proposition}
\begin{proof}
Recall first that any automorphism of $\mathbb{P}^1$ fixes a point. We prove that the same is true for the automorphisms of $\mathbb{P}^1\times\mathbb{P}^1$. Indeed, any such automorphism is of the form
\vspm
 \begin{center}$(u,v) \mapsto (\alpha(u),\beta(v))$ or $(u,v) \mapsto (\alpha(v),\beta(u))$,\end{center} 
 \vspm
 for some $\alpha,\beta \in \Aut(\mathbb{P}^1)=\PGL(2,\K)$. The first automorphism fixes the point $(a,b)$, where $a,b \in \mathbb{P}^1$ are points fixed by respectively $\alpha$ and $\beta$. The second one fixes the point $(c, \beta(c))$, where $c \in \mathbb{P}^1$ is a point fixed by $\alpha\beta$.
 
Blowing-up the fixed point, and blowing-down the strict pull-backs of the two lines of  $\mathbb{P}^1\times \mathbb{P}^1$ passing through the fixed point, we conjugate the automorphism to an automorphism of $\mathbb{P}^2$.\proofend
 \end{proof}

\bigskip

\begin{proposition}
\label{Prp:ConicBundleconjP2}
Any automorphism of finite odd order of some conic bundle (that preserves the c.b. structure) is birationally conjugate to a linear automorphism of $\Pn$.
\end{proposition}
\begin{proof}
Let us denote by $g$ the automorphism of odd order of the conic bundle induced by $\pi:S\rightarrow \mathbb{P}^1$. Recall that the action of $g$ on the fibres of $\pi$ induces an automorphism $\overline{g}$ of $\mathbb{P}^1$ of odd order $m$, whose orbits have all the same size $m$, except for two fixed points. 

Suppose that one fibre $F$ of $\pi$ is singular. The orbit of $F$ by $g$ is thus a set of singular curves $\{F_1,...,F_n\}$ (where $n=1$ or $n=m$). Furthermore, $g$ acts on the set $T$ of irreducible components of the $F_i$'s, whose size is even, equal to $2n$. Since the order of $g$ is odd, the action of $g$ on $T$ has two orbits of size $n$, and two curves of the same orbit do not intersect. This allows us to blow-down one of the two orbits, to obtain a birational $g$-equivariant morphism from the conic bundle to another one, with fewer singular fibres. 

Continuing by this way, we conjugate $g$ to an automorphism of a conic bundle which has no singular fibre.
Since the fibration is smooth, the surface is an Hirzebruch surface $\mathbb{F}_k$, for some integer $k\geq 0$. If $k\geq 1$, choose one fibre $F$ invariant by $g$ (there exist at least two such fibres). Since $F\cong \mathbb{P}^1$, $g$ fixes at least two points of $F$. Blow-up one point of $F$ fixed by $g$ and not lying on the exceptional section of $\mathbb{F}_k$ (the one of self-intersection $-k$); blow-down then the strict pull-back of $F$, to obtain a $g$-equivariant birational map $\mathbb{F}_k\dasharrow \mathbb{F}_{k-1}$. By this way, we may assume that $g$ acts biregularly on $\mathbb{F}_0=\mathbb{P}^1\times\mathbb{P}^1$, and use Proposition \ref{Prp:AutP1P1conjP2} to achieve the proof.\proofend
\end{proof}

\bigskip

\begin{proposition}
\label{Prp:DelPezzo4conjP2}
Any automorphism of finite odd order of a Del Pezzo surface of degree $\geq 4$ is birationally conjugate to a linear automorphism of $\Pn$.
\end{proposition}
\begin{proof} Recall that a Del Pezzo surface is either $\mathbb{P}^1\times\mathbb{P}^1$ or the blow-up of some points in $\mathbb{P}^2$ in general position (i.e.\ such that no irreducible curve of self-intersection $\leq -2$ belongs to the surface).
 
 Suppose that $g$ acts on a Del Pezzo surface $S$ of degree $\geq 4$. By blowing-down some curves (which gives once again a Del Pezzo surface, with a larger degree), we may assume that $g$ acts minimally on $S$. If $S$ is $\mathbb{P}^1\times\mathbb{P}^1$ or $\mathbb{P}^2$, we are done (Proposition \ref{Prp:AutP1P1conjP2}). 
 
 Otherwise, either $g$ preserves a conic bundle structure, or $\rkPic{S}^g=1$ (Proposition~\ref{Prp:Man}). In the first case, $g$ is birationally conjugate to a linear automorphism of $\Pn$ (Proposition \ref{Prp:ConicBundleconjP2}). In the second case, the degree of $S$ divides the order of $g$ (Lemma~\ref{Lem:SizeOrbits}), which is odd by hypothesis.  The only possibilities for the degree of $S$ are thus $5$ or $7$. We study now both cases.
 \begin{itemize}
 \item
 {\it The degree of $S$ is $7$}, i.e.\ $S$ is the blow-up of two distinct points of $\mathbb{P}^2$. In this case, there are three exceptional divisors on $S$, which are the pull-back $E_1$, $E_2$ of the two points, and the strict pull-back of the line of $\mathbb{P}^2$ passing through the two points. This configuration implies that the set $\{E_1,E_2\}$ is invariant by any automorphism of the surface, which is thus birationally conjugate to a linear automorphism of $\Pn$.
  \item
 {\it The degree of $S$ is $5$}, i.e.\ $S$ is the blow-up of four points of $\mathbb{P}^2$, no $3$ being collinear.
We may assume that the points blowed-up are $(1:0:0)$, $(0:1:0)$, $(0:0:1)$ and $(1:1:1)$. The action of the group $\Aut(S)$ of automorphisms  of $S$ on the $5$ sets of $4$ skew exceptional divisors of $S$ gives rise to an isomorphism of $\Aut(S)$ to the group $\Sym_5$. The group $\Aut(S)$ is generated by the lift of the group $\Sym_4$ of automorphisms of $\Pn$ that leaves invariant the four points blowed-up, and by the quadratic transformation $(x:y:z) \dasharrow (yz:xz:xy)$. The nature of $\Aut(S)$ may be found by direct calculation, and is also well-known for many years (see for example \cite{bib:SK}, \cite{bib:Wim}, \cite{bib:JBTh}, \cite{bib:Dol}). Using Lemma~\ref{Lem:SizeOrbits}, the automorphism $g$ with $\rkPic{S}^g=1$ must have order $5$, and is thus birationally conjugate to $(x:y:z) \dasharrow (x(z-y):z(x-y):xz)$, which is birationally conjugate to a linear automorphism of $\Pn$ (see \cite{bib:BeB}).\proofend
 \end{itemize}
\end{proof}

Proposition \ref{Prp:DelPezzo4conjP2} is false for Del Pezzo surfaces of degree at most $3$. We give now some examples:
\begin{example}
\label{Exa:FermatCubic}
Let  $$S_F=\{(w:x:y:z) \in \mathbb{P}^3\ |\ w^3+x^3+y^3+z^3=0\}$$ be the \emph{Fermat cubic surface}, which is a Del Pezzo surface of degree $3$ (see \cite{bib:Kol}, Theorem III.3.5).
The elements 
$$\begin{array}{l}
\rho_1:(w:x:y:z) \mapsto (w:e^{2\im\pi/3}y:z:x)\\
\rho_2:(w:x:y:z) \mapsto (w:e^{4\im\pi/3}y:z:x)\end{array}
$$ 
are automorphisms of $S_F$. For $i=1,2$, the element $(\rho_i)^3:(w:x:y:z) \mapsto (w:e^{i2\im\pi/3}x:e^{i2\im\pi/3}y:e^{i2\im\pi/3}z)$ fixes the elliptic curve which is the intersection of $S$ with the plane $w=0$, and corresponds to an element of order $3$ described in Example \ref{Exa:Order3}. Since $(\rho_i)^3$ is not birationally conjugate to a linear automorphism of $\mathbb{P}^2$, the same occurs for $\rho_i$.
\end{example}
\begin{example}\label{Exa:S15}
Let 
$$S_{15}=\{(w:x:y:z) \in \mathbb{P}(3,1,1,2)\ |\ w^2=z^3+x(x^5+y^5)\},$$
be a special Del Pezzo surface of degree $1$ (see \cite{bib:Kol}, Theorem III.3.5). The element
$$\theta:(w:x:y:z) \mapsto (w:x:e^{2\im\pi/5} y:e^{2\im\pi/3} z)$$
is an automorphism of the surface $S_{15}$ which has order $15$.
Since $\theta^3$ (which is an element described in Example \ref{Exa:Order5}) fixes an elliptic curve, it is not birationally conjugate to a linear automorphism of $\Pn$, and thus the same occurs for $\theta$.
\end{example}

\bigskip

\begin{proposition}
\label{Prp:DelPezzoOdd13}
Let $g$ be some birational map of $\Pn$ of finite odd order $\geq 7$. Then, $g$ is birationally conjugate either to a linear automorphism of $\Pn$, or to one of the elements $\rho_1,\rho_2$ described in Example~\ref{Exa:FermatCubic}, or to one of the elements $\theta,\theta^2,\theta^4,\theta^7,\theta^8,\theta^{11},\theta^{13},\theta^{14}$, where $\theta$ is described in Example~\ref{Exa:S15}.
\end{proposition}
\begin{proof}
As we already mentioned, every birational map of $\Pn$ is birationally conjugate to an automorphism of a rational surface $S$ (see for example \cite{bib:DFE}, Theorem 1.4).
Supposing that the action is minimal (i.e.\ that every $g$-equivariant birational morphism $S\rightarrow S'$ is an isomorphism), either $g$ preserves a conic bundle structure on $S$ or $S$ is a Del Pezzo surface (Proposition \ref{Prp:Man}).\\
In the first case, the automorphism is birationally conjugate to a linear automorphism of $\Pn$, since it has odd order (Proposition \ref{Prp:ConicBundleconjP2}).\\
In the second case, if the surface has degree $\geq 4$, the automorphism is birationally conjugate to a linear automorphism of $\mathbb{P}^2$ (Proposition \ref{Prp:DelPezzo4conjP2}). Otherwise, applying Lemma \ref{Lem:SizeOrbits}, the degree of the surface is $1$ or $3$ and divides the order of the automorphism. We enumerate the possibilities:
\begin{itemize}
\item
{\it the degree of $S$ is $3$, and the order of $g$ is a multiple of $3$}.\\
The linear system $|K_S|$ gives rise to the canonical embedding of $S$ in $\mathbb{P}^3$, whose image is a smooth cubic surface (see for example \cite{bib:Kol}, Theorem III.3.5). Since $g$ leaves invariant the linear system $|K_S|$, it is the restriction of a linear automorphism of $\mathbb{P}^3$. 

Suppose first that $S$ is isomorphic to the Fermat cubic surface $S_F$, whose equation is $w^3+x^3+y^3+z^3=0$, and whose group of automorphisms is $(\Z{3})^3\rtimes \Sym_4$, where $(\Z{3})^3$ is the $3$-torsion of $\PGL(4,\K)$ and $\Sym_4$ is the group of permutations of the variables. Since the order of $g$ is odd and at least $7$, its image in $\Sym_4$ is an element of order $3$. The elements of order $3$ of $\Sym_4$ being all conjugate, $g$ is conjugate to an element of the form
$$(w:x:y:z) \mapsto (w:a y: bz:cx),$$
for some $a,b,c$ in the $3$-torsion of $\K^{*}$. We conjugate $g$ by the automorphism $(w:x:y:z) \mapsto (w:bcx:y:cz)$ of $S_F$ and get the automorphism $(w:x:y:z) \mapsto (w:abc y:z:x)$. Since the order of  $g$ is not $3$, $abc$ is not equal to $1$ and is thus a primitive $3$-th root of unity. The two possible cases give the elements $\rho_1$ and $\rho_2$ described in Example \ref{Exa:FermatCubic}.

We proceed now to the study of general cubic surfaces. We denote by $G$ the group generated by $g$, and by $h$ one of the two elements of order $3$ of $G$. Up to  isomorphism (and to the choice of $h$), three possibilities occur (we use the notation $\omega=e^{2\im\pi/3}$):
\begin{itemize}
\item
{\it The automorphism $h$ is $(w:x:y:z) \mapsto (\omega w:\omega x:y:z)$}.\\
The equation of $S$ is thus $L_3(w,x)+L'_3(y,z)=0$, where $L_3,L_3'$ are homogeneous forms of degree $3$; this implies that the surface is isomorphic to the Fermat cubic surface, a case already studied.
\item
{\it The automorphism $h$ is $(w:x:y:z) \mapsto (\omega w:x:y:z)$}.\\
In this case, $h$ fixes an elliptic curve $\Gamma$, which is the intersection of $S$ with the plane of equation $w=0$ ($h$ corresponds to an element described in Example \ref{Exa:Order3}). Note that $G$ commutes with $h$ so it leaves invariant $\Gamma$ and also the plane $w=0$. The action of the group $G$  on the curve $\Gamma$ must thus be cyclic of odd order at least $3$ and corresponds to the action of a cyclic subgroup of $\PGL(3,\K)$ on a smooth cubic curve. If the action is a translation, it does not have fixed points and corresponds to the action of $(x:y:z) \mapsto (x:\omega y:\omega^2 z)$ on a plane cubic curve of equation $x^3+y^3+w^3+\lambda xyz=0$. But this is not possible, since the group obtained by lifting this action is isomorphic to $(\Z{3})^2$ and thus is not cyclic.
It remains the case of an automorphism of an elliptic curve, which has fixed points. The only possibility is an element of order $3$, that acts on the curve of equation $x^3+y^3+z^3=0$. But this case yields once again the Fermat cubic surface.
\item
{\it The automorphism $h$ is $(w:x:y:z) \mapsto (\omega w:\omega^2x:y:z)$}
\\
We see now that this case is incompatible with the hypothesis on $g$. Writing the equation of $S$ as $F=0$, where $F \in \K[w,x,y,z]$, we have $F(\omega w,\omega^2 x,y,z)=\omega^i F(w,x,y,z)$, for some $i=0,1,2$. The different cases give respectively the equation $F_0$, $F_1$ and $F_2$ listed below (note that since $S$ is non-singular, its equation has degree at least $2$ in each variable, we scale to $1$ the non-zero coefficients by change on variables)
\begin{center}$\begin{array}{lllllll}
F_0&=&(ay^3+by^2z+cyz^2+dz^3)&+&(ey+fz)wx&+&w^3+x^3,\\
F_1&=&(ay^2+byz+cz^2)x&+&(ey+fz)w^2&+&wx^2,\\
F_2&=&(ay^2+byz+cz^2)w&+&(ey+fz)x^2&+&w^2x,\end{array}$\end{center}
where $a,b,c,d,e,f \in \K$. We may suppose that $g$ is diagonal, i.e.\ that $$g(w,x,y,z)=(w:\lambda_x\cdot x:\lambda_y\cdot y:\lambda_z\cdot z),$$
for some $\lambda_x,\lambda_y,\lambda_z \in \K^{*}$.

We prove that the equation of $S$ may not be of the form $F_0$. Suppose that this is the case, and deduce first that $\lambda_x^3=1$. As $S$ is non-singular, the equation is of degree $2$ in each variable so either $a\not=0$ or $b\not=0$ and the same is true for $c$ and $d$.
\begin{itemize}
\item
If $ad\not=0$, then $g$ has order $3$, which is excluded.
\item
If $bc\not=0$, we have $\lambda_y^2\lambda_z=\lambda_y\lambda_z^2=1$, so $g$ is once again of order $3$.
\item
By permuting $y$ and $z$ if necessary, we suppose that $ac\not=0$, which implies that $\lambda_y^3=\lambda_y\lambda_z^2=1$, whence $\lambda_y=\pm \lambda_z$. Since the order of $g$ is odd, this implies that $\lambda_y=\lambda_z$, which implies once again, with the relation $\lambda_y^3=1$ that $g$ has order $3$.
\end{itemize}
We prove now that the equation of $S$ is neither of the form $F_1$. Suppose the converse, and see that the non-singularity of $S$ implies that $ac\not=0$ and either $e$ or $f$ is non-equal to zero. We may assume (up to a permutation of $y$ and $z$) that $ace\not=0$, whence\begin{center}$\lambda_z^2\lambda_x=\lambda_y^2\lambda_x=\lambda_y=\lambda_x^2$.\end{center} These relations imply that $\lambda_x\lambda_y=1$ and $\lambda_x^3=\lambda_y^3=1$. Since $g$ is odd, the relation $\lambda_y^2=\lambda_z^2$ implies $\lambda_z=\lambda_y$, so $g$ is once again of order $3$.

The case of $F_2$ is similar to the case of $F_1$, by exchanging the variables $w$ and $x$.
\end{itemize}
\item
{\it the degree of $S$ is $1$}.\\
The linear system $|-2K_S|$ induces a degree $2$ morphism onto a quadric cone in $Q\subset\mathbb{P}^3$, ramified over the vertex $v$ of $Q$ and a smooth curve $C$ of genus $4$. Moreover $C$ is the intersection of $Q$ with a cubic surface. (See \cite{bib:BaB}, \cite{bib:deF}, \cite{bib:Dol}.)

Note that a quadric cone is isomorphic to the weighted projective plane $\mathbb{P}(1,1,2)$ and the ramification curve $C$ has equation of degree $6$ there. Up to a change of coordinates, we may thus assume that the surface $S$ has the equation
\begin{center}
$w^2=z^3+F_4(x,y)z+F_6(x,y)$
\end{center}
in the weighted projective space $\mathbb{P}(3,1,1,2)$, where $F_4$ and $F_6$ are forms of respective degree $4$ and $6$ (see \cite{bib:Kol}, Theorem III.3.5). Remark that multiple roots of $F_6$ are not roots of $F_4$, since $S$ is non-singular, and the point $v=(1:0:0:1)=(-1:0:0:1)$ is the vertex of the quadric.

The double covering of the quadric $Q \cong \mathbb{P}(1,1,2)$ gives an exact sequence
\begin{center}
$1\rightarrow <\sigma> \rightarrow \Aut(S) \rightarrow \Aut(Q)_{C}$,
\end{center}
where $\Aut(Q)_{C}$ denote the automorphisms of $Q$ that leaves invariant the ramification curve $C=\{(x:y:z)\ |\ z^3+zF_4(x,y)+F_6(x,y)=0\}$. (In fact we can prove that the right homomorphism is surjective, but we will not need it here). A quick calculation shows that any element of $\Aut(Q)_{C}$ belongs to $\mathrm{P}(\GL(2,\K) \times \GL(1,\K))$. This implies that $\Aut(S)\subset \mathrm{P}(\GL(1,\K)\times \GL(2,\K) \times \GL(1,\K))$.

Note that $|K_S|$ is a pencil of elliptic curves, parametrised by the $(x,y)$-coordinates, which has one base point, $v=(1:0:0:1)$. Any automorphism of $S$ acts thus on the elliptic bundle and fixes the vertex $v$ of $Q$. This induces another exact sequence
\begin{center}
$1\rightarrow G_S \rightarrow \Aut(S) \stackrel{\pi}{\rightarrow} \Aut(\mathbb{P}^1)$,
\end{center}
where \begin{center}$G_S=\left\{\begin{array}{llll}<&(w:x:y:z) \mapsto (-w:x:y:e^{2\im\pi/3}z)&>\cong \Z{6} &\mbox{ if }F_4 =0, \\
<&(w:x:y:z) \mapsto (-w:x:y:z)&>\cong\Z{2} & \mbox{ otherwise.}\end{array}\right.$\end{center}
(Note that the involution that belongs to $G_S$ is called "Bertini involution"). Denoting by $G$ the group generated by our automorphism $g$, the group $G \cap G_S$ is either trivial or cyclic of order $3$. We study the two different cases.
\begin{itemize}
\item {\it The group $G \cap G_S$ is trivial.}\\
In this case, the action of $G$ on the elliptic pencil (which is cyclic and diagonal) has the same order as $g$, which is by hypothesis at least $7$. As both $F_4$ and $F_6$ are preserved by this action, both are monomials.

Then, either $y^2$ or $x^2$ divides $F_6$, which implies that $F_4$ is a multiple of $x^4$ or $y^4$. (Recall that the double roots of $F_6$ are not roots of $F_4$). Up to an exchange of coordinates, we may thus suppose that $F_4=x^4$ and $F_6=y^6$ or $F_6=xy^5$.

In the first case, the equation of the surface is $w^2=z^3+x^4z+y^6$, and its group of automorphisms is isomorphic to $\Z{2}\times\Z{12}$, generated by the Bertini involution $(w:x:y:z) \mapsto (-w:x:y:z)$ and $(w:x:y:z) \mapsto (\im w:x:e^{2\im\pi/12}:-z)$. This case is thus not possible, since the order of $g$ is odd and at least $7$.

In the second case, the equation of the surface is $w^2=z^3+x^4z+xy^5$ and its group of automorphisms is isomorphic to $\Z{20}$, generated by $(w:x:y:z) \mapsto (\im w:x:e^{2\im\pi/10}y:-z)$. We obtain once again a contradiction.
\item
{\it The group $G\cap G_S$ is cyclic of order $3$, generated by $(w:x:y:z) \mapsto (w:x:y:e^{2\im\pi/3}z)$}\\
In this case, 
$L_4=0$, so $L_6$ has exactly $6$ distinct roots. Since the action of $G$ on these roots must be of odd order $\geq 3$, it must be of order $3$ or $5$.

The element $g$ is of the form $(w:x:y:z) \mapsto (\lambda_w w:x:\alpha y:\lambda_z z)$, for some $\lambda_w,\lambda_z\in\K^{*}$, where $\alpha$ is a $n$-th root of the unity, and $n=3$ or $n=5$. This implies that $L_6$ is, up to a linear change on $x$ and $y$, respectively $x^6+ax^3y^3+y^6$, for some $a\in \K$,  or $x(x^5+y^5)$. Thus, we have $\lambda_w^2=1$ and $\lambda_z^3=1$. Since the order of $g$ is odd and at least $7$, this shows that $\lambda_w=1$, $n=5$ and $\lambda_z$ is a $3$-th root of the unity. The surface $S$ is thus the surface $S_{15}$ of Example \ref{Exa:S15}, and the automorphism is one power of $\theta$, which has order $15$.\proofend
\end{itemize}
\end{itemize}
\end{proof}

We are now able to prove Theorem \ref{Thm:FiniteOdd}:
\begin{proof}{\it of Theorem \ref{Thm:FiniteOdd}}\\
Using Proposition \ref{Prp:DelPezzoOdd13} above, and the fact that all the linear automorphisms of some given finite order are birationally conjugate (\cite{bib:BaB}), there exists exactly one single conjugacy class of elements of the Cremona group of some given odd order $n\not=3,5,9,15$, which is represented by the linear automorphism $\alpha_n:(x:y:z)\mapsto (x:y:e^{2\im\pi/n}z)$.\\
Using the same results, the elements of order $9$ of the Cremona group are birationally conjugate to one of the three elements $\alpha_9$, $\rho_1$, $\rho_2$, where $\alpha_9$ is the automorphism $\alpha_9:(x:y:z)\mapsto (x:y:e^{2\im\pi/9}z)$ of $\Pn$ and $\rho_1,\rho_2$ are the automorphisms of the Fermat cubic surface $S_F$ described in Example \ref{Exa:FermatCubic}. It remains to show that these three elements are not birationally conjugate. Firstly, since $(\rho_1)^2$ and $(\rho_2)^2$ both fix an elliptic curve, neither of them is birationally conjugate to a linear automorphism of $\Pn$. Thus $\alpha$ is neither conjugate to $\rho_1$, nor to $\rho_2$. Secondly, the elements $(\rho_1)^2$ and $(\rho_2)^2$ are diagonal in $\PGL(4,\K)$ and have distinct eigenvalues (up to multiplication), so are not conjugate by an element of $\PGL(4,\K)$. This implies that $\rho_1,\rho_2$, which are elements of the group of automorphisms of the Fermat cubic surface $S_F$, are not conjugate in this group. Suppose now that these two elements are conjugate by some birational transformation $\varphi$ of $S_F$. Then, since $\varphi$ is $G$-equivariant, where $G\cong \Z{9}$, we may factorise it into a composition of elementary $G$-equivariant links (see \cite{bib:Isk5}). The first link is the blow-up of some orbit of $G$, and conjugates the group generated by $\rho_1$ to a group of automorphisms of some conic bundle. But this is not possible, because any such automorphism is birationally conjugate to a linear automorphism of $\Pn$ (Proposition \ref{Prp:ConicBundleconjP2}). The birational map $\varphi$ is thus an automorphism of the surface, but we proved that $\rho_1$ and $\rho_2$ are not conjugate in $\Aut(S_F)$. These elements are thus neither birationally conjugate. Summing up, there are three conjugacy classes of elements of order $9$ in the Cremona group, represented by $\alpha_{9}$, $\rho_1$ and $\rho_2$.\\
The case of elements of order $15$ is similar. Using once again Proposition \ref{Prp:DelPezzoOdd13} and \cite{bib:BaB}, any element of order $15$ of the Cremona group is birationally conjugate either to $\alpha_{15}:(x:y:z)\mapsto (x:y:e^{2\im\pi/9}z)$, or to one of the generators of the group $<\theta>$, generated by the automorphism $\theta\in \Aut(S_{15})$ described in Example \ref{Exa:S15}. Since the $5$-torsion of $<\theta>$ fixes an elliptic curve, no generator of $<\theta>$ is birationally conjugate to $\alpha_{15}$. Note that the group of automorphisms of $S_{15}$ is isomorphic to $\Z{30}$, generated by $\theta$ and the Bertini involution. Two distinct elements of the group $<\theta> \subset \Aut(S_{15})$ are thus not conjugate by an automorphism of $S_{15}$. The same argument as before shows that the elements are not birationally conjugate. There are thus exactly $9$ conjugacy classes of elements of order $15$ in the Cremona group, represented by $\alpha_{15}$, $\theta$, $\theta^2$, $\theta^4$, $\theta^7$, $\theta^8$, $\theta^{11}$, $\theta^{13}$ and $\theta^{14}$.\proofend
\end{proof}

\end{document}